\documentclass{ifacconf}

\usepackage{graphicx}
\usepackage{natbib}    
\usepackage{appendix}
\usepackage[utf8]{inputenc}
\usepackage{amsmath,amssymb,amsfonts}
\usepackage{mathrsfs}
\usepackage{stmaryrd}
\usepackage{graphicx}
\usepackage{mathtools}
\usepackage{tikz}
\usepackage{enumitem}
\usepackage{array}
\newcommand{\poubelle}[1]{}

\DeclarePairedDelimiter{\abs}{\lvert}{\rvert} 
\DeclarePairedDelimiter{\floor}{\lfloor}{\rfloor} 

\DeclareMathOperator{\diff}{d\!}

\newcommand{\suchthat}{\ifnum\currentgrouptype=16 \mathrel{}\middle|\mathrel{}\else\mid\fi}

\theoremstyle{plain}
\newtheorem{theorem}{Theorem}[section]
\newtheorem{proposition}[theorem]{Proposition}
\newtheorem{corollary}[theorem]{Corollary}

\theoremstyle{definition}
\newtheorem{definition}[theorem]{Definition}
\newtheorem{remark}[theorem]{Remark}

\begin{document}

\setlist[enumerate]{label={(\alph*)}}

\begin{frontmatter}

\title{Pad\'e Approximation and Hypergeometric Functions: A Missing Link with the Spectrum of Delay-Differential Equations} 

\author[First,Second]{Islam Boussaada} 
\author[First]{Guilherme Mazanti}
\author[First]{Silviu-Iulian Niculescu} 
\address[First]{Universit\'e Paris-Saclay, CNRS, CentraleSup\'elec, Inria, Laboratoire des Signaux et Syst\`emes, 91190, Gif-sur-Yvette, France (e-mails: \{Islam.Boussaada, Guilherme.Mazanti, Silviu.Niculescu\}@l2s.centralesupelec.fr).}
\address[Second]{Institut Polytechnique des Sciences Avanc\'{e}es (IPSA), 63 boulevard de Brandebourg, 94200 Ivry-sur-Seine, France}

\begin{abstract}   
It is well known that rational approximation theory involves degenerate hypergeometric functions and, in particular, the Pad\'e approximation of the exponential function is closely related to Kummer hypergeometric functions. Recently, in the context of the study of the exponential stability of the trivial solution of delay-differential equations, a new link between the degenerate hypergeometric function and the zeros distribution of the characteristic function associated with linear delay-differential equations was emphasized. Such a link allowed the characterization of a property of time-delay systems known as \emph{multiplicity-induced-dominancy (MID)}, which opened a new direction in designing low-complexity controllers for time-delay systems by using a partial pole placement idea. Thanks to their relations to hypergeometric functions, we explore in this paper links between the spectrum of delay-differential equations and Pad\'e approximations of the exponential function. This note exploits and further comments recent results from [I.~Boussaada, G.~Mazanti and S-I.~Niculescu. 2022, Comptes Rendus. Math\'ematique] and [I.~Boussaada, G.~Mazanti and S-I.~Niculescu. 2022, Bulletin des Sciences Math\'ematiques].
\end{abstract}


\end{frontmatter}

\section{Introduction}

On beyond of the arithmetic theory and the number theory, the theory of continued fractions is strongly related to the theory of analytic functions; it is involved in topics such as the study of definite integrals, power series and the summation of divergent series, see for instance \cite{Perron,wall2018}. Continuous fractions theory goes back to the pioneering work of T.~S.~Stieltjes mainly elaborated at the end of the 19th century. It appears that several links exist between continuous fractions theory with  concepts and  methods such as the well known Pad\'e rational approximation and the $\epsilon$-algorithm. As a matter of fact, the $\epsilon$-algorithm is a transformation of the partial sums of a series into the corresponding Pad\'e quotients or the corresponding continued fractions. In fact, the $\epsilon$-algorithm  is an old but still relevant way to improve the convergence rate of slowly convergent sequences, see, for instance, \cite{graves} and references therein. Moreover, a strong connection between such an algorithm and the Pad\'e approximation has been underlined by P.~Wynn in the early 60s, \cite{Wynn66}. Through his results on continuous fractions theory, Wynn's work covered several problems in analytic functions theory, such as the distribution of zeros of a class of special functions called the Kummer (degenerate) hypergeometric functions or, equivalently, the \emph{Wittaker hypergeometric functions}, see, e.g., \cite{Wynn1973Zeros}. 

More contemporary contributions have pointed out to a missing link existing between the zeros of such a class of special functions with the spectrum distribution of an appropriate Delay-Differential Equation (DDE) including one pointwise delay, see for instance \cite{MBN-2021-JDE,BMN-2022-BSM}. Interestingly, this link allowed deriving an exhaustive characterization of a property for DDE\footnote{in both retarded as well neutral cases} called \emph{Generic Multiplicity Induced Dominancy} (GMID), see \cite{BMN-2022-CRAS}. Indeed, the GMID is a special case of a more general property called \emph{Multiplicity Induced Dominancy} (MID). Such a property consists of determining the conditions under which a given multiple complex zero of a quasipolynomial is dominant. It should be noted that the GMID property asserts that a DDE spectral value admitting the maximal multiplicity is necessarily the corresponding spectral abscissa. However, multiple roots with intermediate admissible multiplicities may or may not be dominant. Thanks to this property, a consistent control strategy is proposed in \cite{Boussaada2020Multiplicity,Balogh21}, which consists in assigning a root with an admissible multiplicity once appropriate conditions guaranteeing its dominancy have been established. Furthermore, the MID property may be used to tune standard controllers. For instance, in \cite{ma:hal-03427550} it is applied to the systematic tuning of the stabilizing PID controller of a first-order plant. Here, we aim at assigning dominant multiple real roots with admissible codimensions.

Beyond the reminder of the concepts and links discussed, the contribution of the present note is threefold. First, to highlight the existing links between roots of generic maximal multiplicity for a quasipolynomial, Kummer confluent geometric functions, and Pad\'e approximations of the exponential function. Second, the use of Pad\'e approximation allows simplifying some newly obtained analytic functions results, thus opening some new perspectives in a \emph{partial pole placement} methodology for linear time-invariant DDE. Finally, it allows to slightly correct some of the results derived in \cite{Saff-Varga-1978} where the links between Kummer functions and Pad\'e approximations are underlined. 

The remaining of the paper is organized as follows: Section~\ref{Prer} recalls the definitions and properties of some concepts studied in the paper. Next, Section~\ref{MR} is dedicated to the main result. Section~\ref{IE} presents further comments on the distribution of the Kummer and the Whittaker zeros distribution. Some concluding remarks end the paper.

{\bf Notations:} Throughout the paper, the following notations are used:  $\mathbb N^\ast$ denotes the set of positive integers and $\mathbb N = \mathbb N^\ast \cup \{0\}$. The set of all integers is denoted by $\mathbb Z$ and, for $a, b \in \mathbb R$, we denote $\llbracket a, b\rrbracket = [a, b] \cap \mathbb Z$, with the convention that $[a, b] = \emptyset$ if $a > b$. For a complex number $s$, $\Re(s)$ and $\Im(s)$ denote its real and imaginary parts, respectively. The open left and right complex half-planes are the sets $\mathbb C_-$ and $\mathbb C_+$, respectively, defined by $\mathbb C_- = \{s \in \mathbb C \suchthat \Re(s) < 0\}$ and $\mathbb C_+ = \{s \in \mathbb C \suchthat \Re(s) > 0\}$. Given $k, n \in \mathbb N$ with $k \leq n$, the binomial coefficient $\binom{n}{k}$ is defined as $\binom{n}{k} = \frac{n!}{k! (n-k)!}$ and this notation is extended to $k, n \in \mathbb Z$ by setting $\binom{n}{k} = 0$ when $n < 0$, $k < 0$, or $k > n$. For $\alpha \in \mathbb C$ and $k \in \mathbb N$, $(\alpha)_k$ is the \emph{Pochhammer symbol} for the \emph{ascending factorial}, defined inductively as $(\alpha)_0 = 1$ and $(\alpha)_{k+1} = (\alpha+k) (\alpha)_k$.

\section{prerequisites}
\label{Prer}

This section provides a brief presentation of the definitions and results that shall be of use in the sequel concerning Padé rational approximation and degenerate hypergeometric functions.

\subsection{On the Pad\'e rational approximation}
\label{sec:Pade}

The Padé rational approximation can be seen as a generalization of the Taylor expansion of a function in which one wishes to approximate some function $f$ around a given point (which is typically normalized to $0$) by a rational function, i.e., by a function $R$ expressed as $R = \frac{Q}{P}$, where $Q$ and $P$ are polynomials. We start this subsection by a precise definition of what is meant by Padé approximation in this paper.

\begin{definition}
Let $f: U \to \mathbb C$ be an analytic function defined on a neighborhood $U \subset \mathbb C$ of $0$. Given $(m, n) \in \mathbb N^2$, we say that a pair of polynomials $(Q, P)$ with $\deg Q \leq m$ and $\deg P \leq n$ is a \emph{Padé approximation of $f$ of order $(m, n)$} if
\[
P(z) f(z) - Q(z) = \mathscr{O}(z^{n + m + 1}) \qquad \text{ as } z \to 0.
\]
\end{definition}

Classical results on Padé approximations can be found, for instance, in \cite{Baker1996Pade} and \cite[Chapter~4]{Brezinski2002Computational}, and we also refer to \cite{Brezinski1991History} for a historical presentation of this topic. For every $(m, n) \in \mathbb N^2$, the Padé approximation $(Q, P)$ of an analytic function $f$ defined on a neighborhood of $0$ exists and is unique up to a multiplicative constant. If $n = 0$, the corresponding Padé approximation $(Q, 1)$ is such that $Q$ coincides with the Taylor--Maclaurin expansion of $f$ of order $m$.

Padé approximations are widely used in rational approximation theory (see, for instance, \cite{szego1924,dieudonne1935,Saff-Varga-1978,baratchart} and the references therein). In particular, several applications involve Padé approximations of the exponential function $z \mapsto e^z$ and some of the above works highlighted links between such approximations and \emph{degenerate hypergeometric functions}. As reported in \cite[p246]{Perron}, the Padé approximation of order $(m, n)$ of the exponential function is the pair $(Q_{n, m}, P_{n, m})$ defined for $z \in \mathbb C$ by
\begin{equation}\label{PADE}
\left\{
\begin{aligned}
P_{n,m}(z) & = \sum_{k=0}^n\frac{(n+m-k)! n! z^k}{k! m! (n-k)!},\\
Q_{n,m}(z) & = \sum_{k=0}^m\frac{(n+m-k)! (-z)^k}{k! (m-k)!},
\end{aligned}
\right.
\end{equation}
and the rational approximation $R_{n,m}(z)  = \frac{Q_{n,m}(z)}{P_{n,m}(z)}$ satisfies
\begin{equation*}
 R_{n,m}(z)\xrightarrow[n,m\to\infty]{}e^{z},
\end{equation*}
uniformly on compact subsets of $\mathbb C$. Further, as can be found in \cite{Perron}, the Pad\'e remainder $e^z-R_{n,m}(z)$ satisfies
\begin{equation}
\label{eq:Pade-remainder}
P_{n,m}(z) \left(e^z - R_{n,m}(z)\right) = \frac{z^{n+m+1}}{m!}\int_0^1 e^{tz} (1-t)^m t^n \diff t,
\end{equation}
which, as detailed in Section~\ref{sec:Whittaker} below, is closely related to Kummer hypergeometric functions (see for instance \cite{Saff-Varga-1978} and the references therein). We also highlight the fact that, from a control theory point of view, Pad\'e approximations of the exponential function have been exploited in approximating delay systems via finite-dimensional systems (see, for instance, \cite{chmt91:ijc} for a deeper discussion on the approximation of the parametric family of first-order stable delay systems).

\begin{remark}
The expressions of $P_{n, m}$ and $Q_{n, m}$ provided above use a different normalization with respect to the expressions in \cite{Perron}: we have chosen here the normalization consisting in requiring the denominator $P_{n, m}$ to be monic, i.e., the coefficient of the term $z^n$ is $1$.
\end{remark}

\subsection{Degenerate hypergeometric functions}
\label{sec:Whittaker}

We present in this section the definition and main properties of Kummer and Whittaker hypergeometric functions. The presentation in this section follows that of \cite{BMN-2022-BSM}, and we start by providing the definition and main properties of Kummer hypergeometric functions (for further details, see, for instance, \cite{Buchholz1969Confluent, Erdelyi1981Higher, Olver2010NIST} and references therein). For an historical perspective, the reader is referred to \cite{Kam:37}.

\begin{definition}
Let $a, b \in \mathbb C$ and assume that $b$ is not a nonpositive integer. The \emph{Kummer confluent hypergeometric function} $\Phi(a, b, \cdot): \mathbb C \to \mathbb C$ is the entire function defined for $z \in \mathbb C$ by the series
\begin{equation}
\label{DefiConfluent}
\Phi(a, b, z) = \sum_{k=0}^{\infty} \frac{(a)_k}{(b)_k} \frac{z^k}{k!}.
\end{equation}
\end{definition}

The series in \eqref{DefiConfluent} converges for every $z \in \mathbb C$, and the function $\Phi(a, b, \cdot)$ satisfies the \emph{Kummer differential equation}
\begin{equation}
\label{KummerODE}
z \frac{\partial^2 \Phi}{\partial z^2}(a, b, z) + (b - z) \frac{\partial \Phi}{\partial z}(a, b, z) - a \Phi(a, b, z) = 0.
\end{equation}
The equation \eqref{KummerODE} admits two linearly independent solutions, which sometimes are both called Kummer confluent hypergeometric functions. In the present paper, we are concerned only with the solution given by \eqref{DefiConfluent}.

Kummer functions admit the integral representation
\begin{equation}
\label{eq:Kummer-integral}
\Phi(a, b, z) = \frac{\Gamma(b)}{\Gamma(a) \Gamma(b - a)} \int_0^1 e^{zt} t^{a-1} (1-t)^{b-a-1} \diff t,
\end{equation}
for every $a, b, z \in \mathbb C$ such that $\Re(b) > \Re(a) > 0$, where $\Gamma$ denotes the Gamma function. For all complex numbers $a, b, z$ such that $b$ is not a nonpositive integer, Kummer functions also satisfy the relation
\begin{equation}
\label{eq:Kummer-relation}
\Phi(a, b, z) = e^z \Phi(b - a, b, -z),
\end{equation}
and we have in addition
\begin{equation}
\label{eq:Kummer-in-zero}
\Phi(a, b, z) = 1 + \mathscr{O}(z) \quad \text{ as } z \to 0.
\end{equation}

Kummer confluent hypergeometric functions have close links with Whittaker functions, defined as follows (see, e.g., \cite{Olver2010NIST}).
\begin{definition}
Let $k, l \in \mathbb C$ and assume that $2l$ is not a negative integer. The \emph{Whittaker function} $\mathcal M_{k, l}$ is the function defined for $z \in \mathbb C$ by
\begin{equation}
\label{KummerWhittaker}
\mathcal{M}_{k,l}(z) = e^{-\frac{z}{2}} z^{\tfrac{1}{2} + l} \Phi(\tfrac{1}{2} + l - k, 1 + 2 l, z).
\end{equation}
\end{definition}

If $\frac{1}{2} + l$ is not an integer, the function $\mathcal M_{k, l}$ is a multi-valued complex function with branch point at $z = 0$. Whenever $2l$ is not a negative integer, the nontrival roots of $\mathcal M_{k, l}$ coincide with those of $\Phi(\tfrac{1}{2} + l - k, 1 + 2 l, \cdot)$ and $\mathcal M_{k, l}$ satisfies the \emph{Whittaker differential equation}
\begin{equation}
\label{Whittaker}
\varphi''(z) =\left(\frac{1}{4}-\frac{k}{z}+\frac{l^2-\frac{1}{4}}{z^2}\right)\varphi(z).
\end{equation}
Similarly to the Kummer differential equation \eqref{KummerODE}, other solutions of the Whittaker differential equation \eqref{Whittaker} are also known as Whittaker functions in other works, but they will not be used in this paper. Notice also that, since $\mathcal M_{k, l}$ is a nontrivial solution of the second-order linear differential equation \eqref{Whittaker}, any nontrivial root of $\mathcal M_{k, l}$ is necessarily simple. The roots of Whittaker functions satisfy the following immediate symmetry property.

\begin{proposition}
\label{PropWhittakerSymmetricRoots}
Let $k, l \in \mathbb C$ and assume that $2l$ is not a negative integer. If $z \in \mathbb C \setminus \{0\}$ is a nontrivial root of $\mathcal M_{k, l}$, then $-z$ is a root of $\mathcal M_{-k, l}$.
\end{proposition}

The following result from \cite{BMN-2022-BSM}, whose proof relies on the Green--Hille transform of \eqref{Whittaker} (see, e.g., \cite{hille1922}), provides information on the location of nontrivial zeros of Whittaker functions.

\begin{proposition}
\label{PropTsvetkoffCorrected}
Let $k,l\,\in \mathbb R$ be such that $2\,l-1 \geq 0$.
\begin{enumerate}
\item\label{PropTsvetkoffCorrected-Imaginary} If $k=0$, then all nontrivial roots $z$ of $\mathcal M_{k,l}$  are purely imaginary.
\item\label{PropTsvetkoffCorrected-k-geq-0} If $k>0$, then all nontrivial roots $z$ of $\mathcal M_{k,l}$ satisfy $\Re(z) > 0$.
\item\label{PropTsvetkoffCorrected-k-leq-0} If $k<0$, then all nontrivial roots $z$ of $\mathcal M_{k,l}$ satisfy $\Re(z) < 0$.
\item\label{PropTsvetkoffCorrected-k-neq-0} If $k\neq0$, then all nontrivial roots $z$ of $\mathcal M_{k,l}$ satisfy
\begin{equation}
\label{eq:bound}
4 k^2 {\Im(z)}^2 - \left(4({l}^{2}-{k}^{2}) - 1\right) {\Re(z)}^{2} > 0.
\end{equation}
\end{enumerate}
Moreover, in all cases, all non-real roots $z$ of $\mathcal M_{k, l}$ satisfy $\abs{z} > \sqrt{4\,l^2-1}$.
\end{proposition}

As an immediate consequence of \eqref{KummerWhittaker} and Proposition~\ref{PropTsvetkoffCorrected}, we have the following result on the location of zeros of Kummer functions with real parameters, which was stated in \cite{BMN-2022-BSM}.

\begin{corollary}
\label{CorZerosKummer}
Let $a,\,b\in \mathbb R$ be such that $b\geq 2$.
\begin{enumerate}
\item If $b = 2a$, then all nontrivial roots $z$ of $\Phi(a,b, \cdot)$ are purely imaginary. 
\item\label{CorZerosKummer-k-geq-0} If $b > 2a$, then all nontrivial roots $z$ of $\Phi(a,b, \cdot)$ satisfy $\Re(z) > 0$.
\item\label{CorZerosKummer-k-leq-0} If $b < 2a$, then all nontrivial roots $z$ of $\Phi(a,b, \cdot)$ satisfy $\Re(z) < 0$.
\item If $b \neq 2a$, then all nontrivial roots $z$ of $\Phi(a,b, \cdot)$ satisfy $(b - 2a)^2 {\Im(z)}^2 - \left(4a(b-a) - 2b\right) {\Re(z)}^{2} > 0$.
\end{enumerate}
Moreover, in all cases, all non-real roots $z$ of $\Phi(a, b, \cdot)$ satisfy $\abs{z} > \sqrt{b(b - 2)}$.
\end{corollary}

\subsection{Delay-differential equations and their characteristic equations}

Consider the linear time-invariant delay-differential equation
\begin{equation}
\label{eq:DDE}
y^{(n)}(t) + \sum_{k=0}^{n-1} a_k y^{(k)}(t) + \sum_{k=0}^m \alpha_k y^{(k)}(t - \tau) = 0,
\end{equation}
where $y(\cdot)$ is the real-valued unknown function, $\tau > 0$ is the delay, and $a_0, \dotsc, a_{n-1}, \alpha_0, \dotsc, \alpha_{m}$ are real coefficients. Such an equation is said to be of \emph{retarded type} if $m < n$ (i.e., if the highest order of derivation appears only in the non-delayed term $y^{(n)}(t)$), or of \emph{neutral type} if $m = n$. We refer to \cite{Hale1993Introduction, Michiels2014Stability} for more information on delay-differential equations and for proofs of the major facts stated in the sequel.

Characterizing the stability of the zero solution of \eqref{eq:DDE} in terms of the coefficients $a_0, \dotsc, a_{n-1}, \alpha_0, \dotsc, \alpha_{m}$ is in general a challenging question. Exponential stability can be characterized through the roots of the \emph{characteristic function} of \eqref{eq:DDE}, which is the entire function $\Delta: \mathbb C \to \mathbb C$ defined for $s \in \mathbb C$ by
\begin{equation}
\label{eq:Delta}
\Delta(s) = s^n + \sum_{k=0}^{n-1} a_k s^k + e^{-s\tau} \sum_{k=0}^m \alpha_k s^k.
\end{equation}
More precisely, the zero solution of \eqref{eq:DDE} is exponentially stable if and only if there exists $\gamma < 0$ such that, for every root $s$ of $\Delta$, we have $\Re(s) \leq \gamma$. Rightmost roots of $\Delta$ thus determine the asymptotic behavior of the system and its rate of exponential stability.

It should be mentioned that the function $\Delta$ from \eqref{eq:Delta} is a particular case of a \emph{quasipolynomial}. Quasipolynomials have been extensively studied in the literature (see, e.g., \cite{Bellman1963Differential}, \cite{Hale1993Introduction}) and, in particular, a result presented in \cite[Part Three, Problem~206.2]{Polya1998Problems} implies that any root $s$ of $\Delta$ has multiplicity at most $\mathscr{D}_{PS} = n + m + 1$. The number $\mathscr{D}_{PS}$ is known as the \emph{degree} of $\Delta$, and such a bound on the multiplicity of a root of $\Delta$ is known as \emph{P\'olya--Szeg\H{o} bound}.

\section{Main results}\label{MR}

The main result of this paper is the following, which highlights the links between roots of maximal multiplicities of $\Delta$, Kummer confluent hypergeometric functions, and Padé approximations of the exponential function. More precisely, our main result characterizes the case where $\Delta$ admits a root of maximal multiplicity $\mathscr{D}_{PS}$ through a suitable Kummer hypergeometric function by using the results on the Padé approximations of the exponential function recalled in Section~\ref{sec:Pade}.

\begin{theorem}
\label{thm:main}
Consider the quasipolynomial $\Delta$ given by \eqref{eq:Delta} and let $s_0 \in \mathbb R$. Let $P$ and $Q$ be the polynomials defined by
\[
P(s) = s^n + \sum_{k=0}^{n-1} a_k s^k, \quad Q(s) = \sum_{k=0}^m \alpha_k s^k
\]
and define the quasipolynomial $\widetilde\Delta$ and the polynomials $\widetilde P$ and $\widetilde Q$ by setting, for $z \in \mathbb C$,
\begin{align}
\widetilde\Delta(z) & = \tau^n \Delta\left(s_0 + \tfrac{z}{\tau}\right), \label{eq:DeltaDeltaTilde}\\
\widetilde P(z) & = \tau^n P\left(s_0 + \tfrac{z}{\tau}\right), \\
\widetilde Q(z) & = e^{-s_0 \tau} \tau^n Q\left(s_0 + \tfrac{z}{\tau}\right)
\end{align}
Then the following assertions are equivalent.
\begin{enumerate}
\item\label{item:max-mult} The real number $s_0$ is a root of multiplicity $\mathscr{D}_{PS} = n + m + 1$ of $\Delta$.
\item\label{item:Pade} The pair of polynomials $(-\widetilde Q, \widetilde P)$ forms a Padé approximation of order $(m, n)$ of the exponential function $z \mapsto e^z$.
\item\label{item:DeltaTilde} The quasipolynomial $\widetilde\Delta$ is given by
\begin{equation}
\label{eq:DeltaTilde-as-Kummer}
\widetilde\Delta(z) = \frac{n! z^{n + m + 1}}{(n + m + 1)!} \Phi(m + 1, n + m + 2, -z).
\end{equation}
\item\label{item:Coeffs} The coefficients $a_0, \dotsc, a_{n-1}, \alpha_0, \dotsc, \alpha_m$ of $\Delta$ are given by
\begin{equation}
\label{Coeffs}
\left\{
\begin{tabular}{@{} r @{} l @{}}
$a_k$ & $\displaystyle {} = \left( -1 \right) ^{n-k}n!\,\sum _{j=k}^{n}{\frac {\binom{j}{k} \binom{m+n-j}{m} s_0^{j-k}}{j!\,{\tau}^{n-j}}}$ \\
& \hspace*{\fill} \textnormal{for} $k \in \llbracket 0, n-1\rrbracket$, \\
$\alpha_k$ & $\displaystyle {} = \left( 
-1 \right) ^{n-1}{{e}^{s_{{0}}\tau}}\sum _{j=k}^{m}{\frac {
 \left( -1 \right) ^{j-k} \left( m+n-j \right) !\,s_0^{j-k
}}{k!\, \left( j-k \right) !\, \left( m-j \right) !\,{\tau}^{n-j}}}$ \\
& \hspace*{\fill} \textnormal{for} $k \in \llbracket 0, m\rrbracket$.
\end{tabular}
\right.
\end{equation}
\end{enumerate}
In addition, if the above equivalent assertions are satisfied and $m \leq n$, then $s_0$ is the rightmost root of $\Delta$, i.e., we have $\Re(s) \leq \Re(s_0)$ for every $s \in \mathbb C$ such that $\Delta(s) = 0$.
\end{theorem}

\begin{pf}
Clearly, thanks to the definition of $\widetilde\Delta$, assertion \ref{item:max-mult} is equivalent to requiring $0$ to be a root of multiplicity $\mathscr{D}_{PS}$ of $\widetilde\Delta$. Since $\widetilde\Delta$ is an analytic function and $\mathscr{D}_{PS}$ is the maximal possible multiplicity of any of its roots, the latter fact is equivalent to requiring that $\widetilde\Delta(z) = \mathscr{O}(z^{n + m + 1})$ as $z \to 0$. Using that $\widetilde\Delta(z) = \widetilde P(z) + e^{-z} \widetilde Q(z)$, the latter is equivalent to
\[
\widetilde P(z) e^{z} + \widetilde Q(z) = \mathscr{O}(z^{n + m + 1}) \qquad \text{ as } z \to 0,
\]
which, thanks to the definition of Padé approximation, is equivalent to \ref{item:Pade}.

Using the uniqueness of the Padé approximation of order $(m, n)$ with a monic denominator, we deduce that \ref{item:Pade} is equivalent to requiring that $-\widetilde Q = Q_{n, m}$ and $\widetilde P = P_{n, m}$, where $P_{n ,m}$ and $Q_{n, m}$ are defined in \eqref{PADE}, i.e., it is equivalent to requiring that, for every $z \in \mathbb C$,
\[
\widetilde\Delta(z) = P_{n, m}(z) - e^{-z} Q_{n, m}(z).
\]
In particular, equivalence with \ref{item:Coeffs} follows after straightforward but long computations by using \eqref{eq:DeltaDeltaTilde} (see, e.g., \cite[Lemma~4.2]{MBN-2021-JDE} for detailed computations in the case $m = n-1$, which can be easily generalized to any $m$).

Finally, if \ref{item:Pade} is satisfied, then, combining \eqref{eq:Pade-remainder}, \eqref{eq:Kummer-integral}, and \eqref{eq:Kummer-relation}, we get that
\begin{align*}
e^z \widetilde\Delta(z) & = P_{n, m}(z) (e^z - R_{n, m}(z)) \\
& = \frac{n! z^{n + m + 1}}{(n + m + 1)!} \Phi(n + 1, n + m + 2, z) \\
& = e^z \frac{n! z^{n + m + 1}}{(n + m + 1)!} \Phi(m + 1, n + m + 2, -z),
\end{align*}
yielding \ref{item:DeltaTilde}. Conversely, if \ref{item:DeltaTilde} holds, then, using \eqref{eq:Kummer-in-zero} and \eqref{eq:DeltaTilde-as-Kummer}, we deduce that
\[
\widetilde P(z) e^z + \widetilde Q(z) = e^z \widetilde\Delta(z) = \mathscr{O}(z^{n + m + 1}),
\]
yielding \ref{item:Pade} thanks to the definition of Padé approximation.
\end{pf}

\begin{remark}
The equivalence between items \ref{item:max-mult}, \ref{item:DeltaTilde}, and \ref{item:Coeffs} of Theorem~\ref{thm:main} was already stated and proved in \cite{BMN-2022-CRAS} by exploiting the integral representation \eqref{eq:Kummer-integral}. By exploiting the Padé approximation of the exponential function and in particular \eqref{eq:Pade-remainder}, Theorem~\ref{thm:main} provides the additional equivalence with \ref{item:Pade} and allows for a much simpler proof of the equivalences between \ref{item:max-mult}--\ref{item:Coeffs}.
\end{remark}

Note that, by considering the first equation in \eqref{Coeffs} with $k = n-1$, one obtains the simple and interesting relation between $s_0$, $\tau$, and $a_{n-1}$ given by
\begin{equation}
\label{RelationS0ANMinus1Tau}
s_0 = -\frac{a_{n-1}}{n} - \frac{m+1}{\tau}.
\end{equation}

As remarked in \cite{BMN-2022-CRAS}, if any of the equivalent assertions of Theorem~\ref{thm:main} is satisfied, then $s_0$ is the unique real root of $\Delta$, and, more precisely, it is the unique root of $\Delta$ on the horizontal strip $\{s \in \mathbb C \suchthat \abs{\Im(s)} < \frac{2\pi}{\tau}\}$ of the complex plane. In addition, the assumption in Theorem~\ref{thm:main} of requiring $s_0$ to be real is justified by the fact that nonreal roots of $\Delta$ cannot have a multiplicity equal to the P\'{o}lya--Szeg\H{o} bound $\mathscr D_{PS}$, since any root $s_0$ of $\Delta$ attaining the maximal multiplicity $\mathscr D_{PS}$ necessarily satisfies \eqref{RelationS0ANMinus1Tau}, and thus it will be real since $a_{n-1}$ is real. As an immediate consequence of  Theorem~\ref{thm:main}, if  \eqref{Coeffs} is satisfied for some $s_0 \in \mathbb R$, then the trivial solution of \eqref{eq:DDE} is exponentially stable if and only if $a_{n-1}>-\frac{n\,(m+1)}{\tau}$.

Note that Theorem~\ref{thm:main} provides necessary and sufficient conditions for a real number $s_0$ to be a root of maximal multiplicity of the quasipolynomial $\Delta$ from \eqref{eq:Delta}. The main result of this section states that, under those conditions, $s_0$ is necessarily a dominant root of $\Delta$.

\section{Additional comments on the location of roots of Whittaker functions}\label{IE}

Proposition~\ref{PropTsvetkoffCorrected} is actually a corrected version of a result by G.~E.~Tsvetkov, \cite[Theorem~7]{Tsvetkov1}. In \cite{Tsvetkov1}, the author provides statements of results concerning the location of zeros of Whittaker functions without the corresponding proofs, pointing only to the technique by E.~Hille from \cite{hille1922}. It turns out, as highlighted in \cite[Counterexample~3.1]{BMN-2022-BSM}, that there are counterexamples to the statement of \cite[Theorem~7]{Tsvetkov1}. Proposition~\ref{PropTsvetkoffCorrected}, proved in \cite{BMN-2022-BSM}, uses the technique suggested by G.~E.~Tsvetkov in his paper, based on the Green--Hille transform, to provide a corrected version of \cite[Theorem~7]{Tsvetkov1}.

Some works in the literature have relied on \cite[Theorem~7]{Tsvetkov1} to provide additional properties on the location of roots of Whittaker functions. This is the case, for instance, of \cite[Proposition~3.2]{Saff-Varga-1978}, which studies Whittaker functions motivated by the previously described links with Padé approximations of the exponential function. More precisely, the statement of \cite[Proposition~3.2]{Saff-Varga-1978} is the following.

\begin{proposition}[{\cite{Saff-Varga-1978}}]
\label{PropSaffVarga}
Let $k \in \mathbb R$, $l > 0$, and $z$ be a nontrivial zero of $\mathcal M_{k, l}$.
\begin{enumerate}
\item\label{item:SV-k-geq-0} If $k > 0$, then $\Re(z) > 2k$ and $\Im(z)^2 > \Re(z)^2 \tfrac{4 l^2 - 4k^2 - 1}{4 k^2}$.
\item\label{item:SV-k-leq-0} If $k < 0$, then $\Re(z) < 2k$ and $\Im(z)^2 > \Re(z)^2 \tfrac{4 l^2 - 4k^2 - 1}{4 k^2}$.
\item If $k = 0$, then $\Re(z) = 0$ and $\Im(z)^2 > 4 l^2 - 1$.
\end{enumerate}
\end{proposition}

It turns out that the counterexample to \cite[Theorem~7]{Tsvetkov1} from \cite[Counterexample~3.1]{BMN-2022-BSM} also provides a counterexample to items \ref{item:SV-k-geq-0} and \ref{item:SV-k-leq-0} of Proposition~\ref{PropSaffVarga}. For the sake of completeness, let us detail the counterexample.

Let $l > 0$ and take $k = l + \frac{3}{2}$. If $z$ is a nontrivial root of $\mathcal M_{k, l}$, then, by \eqref{KummerWhittaker}, $z$ is a nontrivial root of $\Phi(-1, 1 + 2 l, \cdot)$. From \eqref{DefiConfluent}, we have $\Phi(-1, 1 + 2 l, z) = 1 - \frac{z}{1 + 2 l}$, and its unique root is $z = 1 + 2l$. In particular, $\Re(z) = 1 + 2 l < 2k = 3 + 2l$, and thus Proposition~\ref{PropSaffVarga}\ref{item:SV-k-geq-0} is not verified. In addition, as a consequence of Proposition~\ref{PropWhittakerSymmetricRoots}, $z = -1 - 2l$ is the unique nontrivial root of $\mathcal M_{-k, l}$ and, since $\Re(z) = -1 - 2l > -2k = -3 - 2l$, Proposition~\ref{PropSaffVarga}\ref{item:SV-k-leq-0} is not verified.

\begin{figure}[ht]
\centering
\resizebox{\columnwidth}{!}{\input{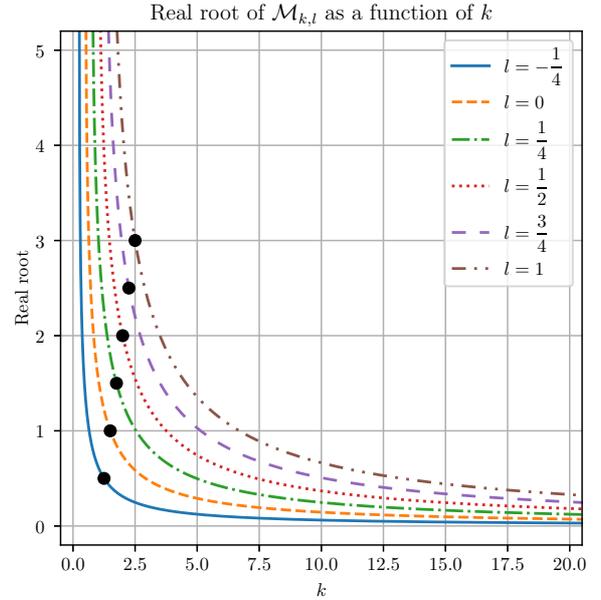}}
\caption{Real root $z(k)$ of $\mathcal M_{k, l}$ satisfying $z(l + \frac{3}{2}) = 1 + 2 l$ for six different values of $l$.}
\label{FigRootWhittaker}
\end{figure}

In view of Proposition~\ref{PropSaffVarga} and Proposition~\ref{PropTsvetkoffCorrected}, a natural question is whether one may provide, in the case $k > 0$, a lower bound on the real part of nontrivial roots of $\mathcal M_{k, l}$ as a function of $k$. A numerical exploration of such a question was provided in \cite{BMN-2022-BSM} by considering the root $z = 1 + 2 l$ of $\mathcal M_{k, l}$ for $k = l + \frac{3}{2}$. Note that such a root exists whenever $l > -\frac{1}{2}$ and it is real and simple. Hence, for every $l > -\frac{1}{2}$, there exists an interval $I_l$ containing $1 + 2 l$ and a curve $k \mapsto z_l(k) \in \mathbb R$ defined on $I_l$ such that, for every $k \in I_l$, $z_l(k)$ is a real root of $\mathcal M_{k, l}$ with $z_l(l + \frac{3}{2}) = 1 + 2l$. The article \cite{BMN-2022-BSM} provides the numerical computation of these curves for $l \in \left\{-\frac{1}{4}, 0, \frac{1}{4}, \frac{1}{2}, \frac{3}{4}, 1\right\}$, which we reproduce in Figure~\ref{FigRootWhittaker}. The black dots $(k, z)$ correspond to $k = l + \frac{3}{2}$ and the root $z = 1 + 2 l$.

An inspection of Figure~\ref{FigRootWhittaker} leads to the conjecture that the maximal interval $I_l$ on which $z_l$ is defined is $I_l = (l + \frac{1}{2}, +\infty)$ and that $z_l(k) \to 0$ as $k \to +\infty$ and $z_l(k) \to +\infty$ as $k \to l + \frac{1}{2}$. In particular, if this conjecture is true, then one cannot expect to correct Proposition~\ref{PropSaffVarga}\ref{item:SV-k-geq-0} by replacing the term $2k$ by any function of $k$ which remains lower bounded as $k \to +\infty$.

{Finally, we notice the remarkable case depicted in \cite{BMN-2022-CRAS}, where one is able to characterize the zeros of the Whittaker function $\mathcal M_{0, n+\frac{1}{2}}$ for $n\in\mathbb{N}$ as $\{i {\zeta} \suchthat \zeta \in \Xi_n\}$, where $\Xi_n$ is the set of $\zeta \in \mathbb R$ satisfying
\begin{equation*}
\tan\left(\frac{\zeta}{2}\right) = \frac{\displaystyle \zeta \sum_{\ell = 0}^{\floor*{\frac{n-1}{2}}} (-1)^\ell \frac{(2n - 2\ell - 1)!}{(2\ell + 1)! (n - 2\ell - 1)!} \zeta^{2\ell}}{\displaystyle \sum_{\ell = 0}^{\floor*{\frac{n}{2}}} (-1)^\ell \frac{(2n - 2\ell)!}{(2\ell)! (n - 2\ell)!} \zeta^{2\ell}}.
\end{equation*} 
The particular case $n=1$, $\Xi_1 = \left\{\zeta \in \mathbb R \suchthat \tan\left(\frac{\zeta}{2}\right) = \frac{\zeta}{2}\right\}$, had already been identified in \cite{MBNC-2021-MTNS}.}

\section{Concluding remarks}

This paper discusses the existing links between roots of generic maximal multiplicity for a quasipolynomial, Kummer confluent geometric functions and Pad\'e approximations of the exponential function. More precisely, we have shown that in the case of a real characteristic root having the maximum multiplicity of a simple quasipolynomial, the corresponding polynomials define a pair of Pad\'e approximations. Moreover, the multiplicity of the root is equal to the P\'olya--Szeg\H{o} bound and it is dominant in the sense that it explicitly defines the spectral abscissa of the spectrum of the corresponding delay-differential equations. As a byproduct of the analysis, we have slightly corrected some of the results from the literature on the Pad\'e approximation.

\bibliography{Bib}

\end{document}